\documentclass{amsart}
\usepackage{verbatim}
\usepackage{amsfonts,amsmath,amsthm,amssymb}
\usepackage[lite]{amsrefs}
\newtheorem{theorem}{Theorem}

\newtheorem{prop}{Proposition}
\newtheorem{cjn}{Conjecture}
\newtheorem{qs}{Question}

\newcommand{\ZZ}{\mathbb{Z}}

\newcommand{\QQ}{\mathbb{Q}}

\def\m1{^{-1}}

\newcommand{\ds}{\displaystyle}

\begin{document}

\title[Torsion Units]
{Zassenhaus Conjecture for Integral Group Ring of Simple Linear Groups}

\author{Joe Gildea}

\address{Department of Mathematics\\
Faculty of Science and Engineering\\
University of Chester\\
England}
\email{j.gildea@chester.ac.uk}

\keywords{Zassenhaus Conjecture, torsion unit, partial augmentation, integral group ring}
\subjclass{16S34, 16U60}

\begin{abstract}
We prove that the Zassenhaus conjecture is true for $PSL(2,8)$ and $PSL(2,17)$. This is a continuation of research initiated by
W. Kimmerle, M. Hertweck and C. H{\"o}fert.
\end{abstract}

\maketitle

\section{Introduction}
Let $U(\ZZ G)$  be the group of units of the integral group ring $\ZZ G$ of a finite group $G$. Then
\[
U(\ZZ G)=\{\pm 1\}\times V(\ZZ G),
\]
where $V(\ZZ G)$ is the group of units of augmentation one.  Throughout this paper, $G$ is always a finite group, torsion units will always represent
torsion units in $ V(\ZZ G)\setminus \{1\}$. One of the most important conjecture in the theory of integral group rings is

\begin{cjn} If $G$ is a finite group, then for each  torsion unit $u\in V(\ZZ G)$ there exists $g\in G$, such that $|u|=|g|$.
\end{cjn}

Hans Zassenhaus formulated a stronger version of this conjecture in \cite{HZ}, which states that:

\begin{cjn} A torsion unit in $V(\ZZ G)$ is said to be rationally conjugate to a group element if it is
conjugate to an element of $G$ by a unit of the rational group ring $\QQ G$.
\end{cjn}

This conjecture was confirmed for Nilpotent Groups in \cite{W,RS} and also some classes of solvable groups (see \cite{MH4} for further details). However these techniques do not transfer to simple groups. The main investigative tool for simple group in relation to the Zassenhaus conjecture is the Luthar-Passi Method (which was introduced in \cite{LP1}). It was confirmed true for all groups up to order $71$ in \cite{HK}. Also this conjecture was validated for $A_5$, $S_5$, central extensions of $S\sb 5$ and other simple finite groups in \cite{LP1,LP2,BH,BHK}. In \cite{MS} partial results were given for $A_6$, and the remaining cases were dealt with in \cite{MH2}. Additionally it was proved for $PSL(2,p)$ when $p=\{7,11,13\}$ in \cite{MH3} and further results regarding $PSL(2,p)$ can be found in \cite{HHK}.

Let $H$ be a  group with  torsion part $t(H)$ of finite exponent.  Let  $\# H$ be the set of primes dividing the order of elements from $t(H)$. The prime graph
of $H$ (denoted by $\pi(H)$) is a graph with vertices labeled by primes from $\# H$, such that vertices $p$ and $q$ are adjacent if and only if there is an element of order $pq$ in the group. The following question which was composed as a problem in \cite{OB} (Problem $37$):

\begin{qs}(Prime Graph Question) If $G$ is a finite group, then $\pi(G)=\pi(V(\ZZ G))$.
\end{qs}

It was established in \cite{WK} that this question is upheld for Frobenius and Solvable groups. It was also confirmed for some Sporadic Simple groups in \cite{BKS,BKM11, BJK,BKL, BKM,BKHS,BKR,BGK,BKM23,BKM24,BKLCW,BKML}. The prime graph question is an intermediate step towards the Zassenhaus conjecture.

We combine the Luthar-Passi Method together with techniques developed in \cite{MH1,MH3} to obtain our results. Our main results are the following:

\begin{theorem}
Let $G= PSL(2,8)$. If $u$ be a torsion unit of $V(\ZZ G)$, then $u$ is rationally conjugated to some element $g \in G$. Hence the Zassenhaus Conjecture is true for $G$.
\end{theorem}

\begin{theorem}
Let $G= PSL(2,17)$. If $u$ be a torsion unit of $V(\ZZ G)$, then $u$ is rationally conjugated to some element $g \in G$. Hence the Zassenhaus Conjecture is true for $G$.
\end{theorem}

Let $u=\sum a_g g$ be a torsion unit of $V(\ZZ G)$, then the order of $u$ divides the exponent of the group (see \cite{CL}). For an element $x \in G$, let $\sum_{g \in X^G} a_g$ be the partial augmentation (denoted by $\varepsilon_C(u)$) of $u$ with respect to it's conjgacy classes in $G$. Let $\nu_i= \varepsilon_{C_i}(u)$ be the $i$-th partial augmentation of $u$.  G. Higman and S. D. Berman \cite{AB} proved that $\nu_1=0$ and $\nu_j=0$ if the conjugacy class $C_j$ consists of a central elements.  Therefore for any torsion $u$ (since $u \in V(\ZZ G)$), we can see that
\[ \nu_2+\nu_3+\cdots+\nu_l=1\]

\noindent where $l$ denotes the class number of $G$.

\begin{prop}\label{CL}(\cite{CL})
Let $u$ be  torsion unit of $V(\ZZ G)$.  The order of $u$ divides the exponent of $G$.
\end{prop}

The following Propositions provide relationships between the partial augmentations and the order
of a torsion unit.

\begin{prop}\label{P:1}(Proposition 3.1 in \cite{MH1} ) Let $u$ be a torsion unit of $V(\ZZ G)$. Let $C$ be a conjucay class of $G$. If $p$ is a prime dividing the order of a representative of $C$ but not the order of $u$ then the partial augmentation $\varepsilon_C(u)=0$.

\end{prop}

\begin{prop}\label{P:2}(Proposition 2.2 in \cite{MH3})
Let $G$ be a finite group and let $u$ be a torsion unit in $V(\ZZ G)$.
\begin{itemize}
\item[$(i)$] If $u$ has order $p^n$, then $\varepsilon_x(u)=0$ for every $x$ of $G$ whose $p$-part is of order
strictly greater than $p^n$.
\item[$(ii)$] If $x$ is an element of $G$ whose $p$-part, for some prime, has order strictly greater than the order of the $p$-part of $u$, then $\varepsilon_x(u)=0$.
\end{itemize}
\end{prop}

\begin{prop}\label{P:3}(\cite{LP1} and Theorem 2.5 in \cite{MRSW})
Let $u$ be a torsion unit of $V(\ZZ G)$ of order $k$. Then $u$ is conjugate in $\QQ G$ to an element $g \in G$ iff
for each $d$ dividing $k$ there is precisely one conjugacy class $C_{i_d}$ with partial augmentation
$\varepsilon_{C_{i_d}}(u^d) \neq 0$.
\end{prop}

\begin{prop}\label{6:1}(Proposition 6.1 in \cite{MH3})
Let $G = PSL(2,p^f )$ for an odd prime $p$, and $f \leq 2$. Then units of order $p$ in $V(\ZZ G)$ are rationally conjugate to elements of $G$.
\end{prop}
\begin{prop}\label{6:3}(Proposition 6.3 in \cite{MH3})
Let $G = PSL(2,p)$ for an odd prime $p$. If the order of a torsion unit $u$ in $V(\ZZ G)$ is divisible by $p$, then $u$ is of order $p$.
\end{prop}
\begin{prop}\label{6:4}(Proposition 6.4 in \cite{MH3})
Let $G = PSL(2,p^f )$, and let $r$ be a prime different from $p$. Then any torsion unit $u$ in $V(\ZZ G)$ of order $r$ is rationally conjugate to an element of $G$.
\end{prop}
\begin{prop}\label{6:6}(Proposition 6.6 in \cite{MH3})
Let $G = PSL(2,p^f)$ with $p \neq  2, 3$. Then any torsion unit in $V(\ZZ G)$ of order $6$ is rationally conjugate to an element of $G$.
\end{prop}

\begin{prop}\label{6:7}(Proposition 6.7 in \cite{MH3})
Let $G = PSL(2, p^f )$. Then for any $p$-regular torsion unit $u$ in
$V(\ZZ G)$, there is an element of $G$ of the same order as $u$.
\end{prop}

For any character $\chi$ of $G$ and any torsion unit of $V(\ZZ G)$ , clearly $\chi(u)=\sum_{i=2}^{l}\nu_i\chi(h_i)$
where $h_i$ is a representative of a conjugacy class $C_i$.


\begin{prop}\label{P:4}(Theorem 1 in \cite{LP1} and \cite{MH3})
Let $p$ be equal to zero or a prime divisor of $|G|$. Suppose that $u$ is an element of $V(\ZZ G)$ of order $k$. Let $z$ be a primitive $k-th$ root of unity. Then for every integer $l$ and any character $\chi$ of $G$, the number
\[ \mu_l(u,\chi,p)=\frac{1}{k} \sum_{d \mid k} Tr_{\QQ(z^d)/\QQ}\{\chi(u^d)z^{-dl}  \}\]
is a non-negative integer.
\end{prop}
We will use the notation $\mu_l(u,\chi,*)$ when $p=0$. The LAGUNA package \cite{lag} for the GAP system \cite{GAP} is a very useful tool when calculating $\mu_l(u,\chi,p)$.

\section{Proof of Theorem 1}

Let $G=PSL(2,8)$. Clearly $|G|=504 = 2^{3} \cdot 3^{2} \cdot 7$ and $exp(G) = 126 = 2 \cdot 3^{2} \cdot 7$. Initially for any torsion unit
of $V(\ZZ G)$ of order $k$, we have that
 \[ \nu_{2a}+ \nu_{3a}+ \nu_{7a}+ \nu_{7b}+ \nu_{7c}+ \nu_{9a}+\nu_{9b}+ \nu_{9c}=1.\]

By Proposition \ref{CL}, we need only to consider torsion units of $V(\ZZ G)$ of order $2$, $3$, $7$, $9$, $6$, $14$, $21$.\\

Case $(i)$. Let $u$ be a torsion unit in $V(\ZZ G)$ of order $2$. By Proposition \ref{P:1}, $\nu_{kx}=0 \; \forall \;kx \in \{3a,7a,7b,7c,9a,9b,9c\}$. Therefore $u$ is rationally conjugated to some element $g \in G$ by Proposition \ref{P:3}.\\

Case $(ii)$. Let $u$ be a torsion unit in $V(\ZZ G)$ of order $3$ and $7$.  Clearly $u$ is rationally conjugated to some element $g \in G$ by Proposition \ref{6:4}.\\

\noindent It remains only to consider torsion units of order $6$, $9$, $14$ and $21$.\\

Case $(iii)$. Let $u$ be a torsion unit in $V(\ZZ G)$ of order $6$.  Using Propositions \ref{P:1} \& \ref{P:2},
$\nu_{2a}+\nu_{3a}=1$. Now $\chi(u^{3}) = \chi(2a)$ and $\chi(u^{2}) = \chi(3a)$.  Applying Proposition \ref{P:4} to the character table (Table 1) , we obtain

\[
\begin{split}
\mu_{3}(u,\chi_{2},*) & = \textstyle \frac{1}{6} (2\gamma_1 +4) \geq 0; \quad
\mu_{0}(u,\chi_{2},*)  = \textstyle \frac{1}{6} (-2\gamma_1 +2) \geq 0; \\
\mu_{1}(u,\chi_{2},*) & = \textstyle \frac{1}{6} (-\gamma_1 +10) \geq 0; \quad
\mu_{0}(u,\chi_{3},*)  = \textstyle \frac{1}{6} (-2 \gamma_2 +8) \geq 0,
\end{split}
\]

\noindent where $\gamma_1 =   \nu_{2a} + 2 \nu_{3a} $ and $\gamma_2= \nu_{2a} -  \nu_{3a} $.  Clearly $\gamma_1=-2$. It follows that there are no
possible integer solutions for $(\nu_{2a},\nu_{3a})$.\\

Case $(iv)$. Let $u$ be a torsion unit in $V(\ZZ G)$ of order $9$.  Using Propositions \ref{P:1} \& \ref{P:2},
$\nu_{3a}+\nu_{9a}+\nu_{9b}+\nu_{9c}=1$. Now $\chi(u^{3}) = \chi(3a)$.  Applying Proposition \ref{P:4} to the character tables (Tables 1, 2 \& 3) , we obtain

\[
\begin{split}
\mu_{2}(u,\chi_{2},7) & =\textstyle \frac{1}{9} (3\gamma_1 + 3)  \geq 0; \quad
\mu_{0}(u,\chi_{2},7)  = \textstyle \frac{1}{9} (-6\gamma_1+ 3) \geq 0; \\
\mu_{1}(u,\chi_{2},2) & = \textstyle \frac{1}{9} (3\gamma_2 + 3) \geq 0; \quad
\mu_{1}(u,\chi_{3},*)  = \textstyle \frac{1}{9} (- 3\gamma_2 + 6) \geq 0;\\
\mu_{2}(u,\chi_{3},*) & = \textstyle \frac{1}{9} ( 3\gamma_3 + 6) \geq 0; \quad
\mu_{2}(u,\chi_{2},2)  = \textstyle \frac{1}{9} (- 3\gamma_3 + 3) \geq 0; \\
\mu_{4}(u,\chi_{2},2) & = \textstyle \frac{1}{9} (- 3 \gamma_4 + 3) \geq 0,
\end{split}
\]

\noindent where $\gamma_1 =  2 \nu_{3a} - \nu_{9a} - \nu_{9b} - \nu_{9c} $, $\gamma_2=2 \nu_{9a} - \nu_{9b} - \nu_{9c}$, $\gamma_3= \nu_{9a} - 2 \nu_{9b} + \nu_{9c}$ and $\gamma_4=\nu_{9a} +\nu_{9b} -2 \nu_{9c}$.  Clearly $\gamma_1 = -1$, $\gamma_2 \in \{ -1, 2 \}$ and $\gamma_3 \in \{ -2, 1\}$.  It follows that
the only possible integer values for $(\nu_{3a},\nu_{9a},\nu_{9b},\nu_{9c})$ are $(0,1,0,0)$, $(0,0,1,0)$ and $(0,0,0,1)$. Therefore $u$ is rationally conjugated to some element $g \in G$ by Proposition \ref{P:3}.\\

Case $(v)$. Let $u$ be a torsion unit in $V(\ZZ G)$ of order $14$.  Using Propositions \ref{P:1} \& \ref{P:2},
$\nu_{2a}+\nu_{7a}+\nu_{7b}+\nu_{7c}=1$. Now $\chi(u^{7}) = \chi(2a)$ and $\chi(u^{2}) = \chi(7a)$.  Applying Proposition \ref{P:4} to the character table (Table 1) , we obtain

\[
\begin{split}
\mu_{7}(u,\chi_{2},*) & = \textstyle \frac{1}{14} (6 \nu_{2a} +8) \geq 0;  \quad   \mu_{0}(u,\chi_{2},*) = \textstyle \frac{1}{14} (-6 \nu_{2a} +6) \geq 0;\\
\mu_{1}(u,\chi_{2},*) & = \textstyle \frac{1}{14} (- \nu_{2a} +8) \geq 0.
\end{split}
\]

\noindent Clearly $\nu_{2a} = 1$. It follows that there are no
possible integer solutions for $\nu_{2a}$. The exact same inequalities are obtained when
$\chi(u^{7}) = \chi(2a)$ and $\chi(u^{2}) = \chi(7b)$, and also when $\chi(u^{7}) = \chi(2a)$ and $\chi(u^{2}) = \chi(7c)$.\\

Case $(vi)$. Let $u$ be a torsion unit in $V(\ZZ G)$ of order $21$.  Clearly $u$ is rationally conjugated to some element $g \in G$ by Proposition \ref{6:7}. This completes the proof.

\vspace{-0.15in}

\begin{table}[h!]
  \begin{center}
 \caption{Character Table of $PSL(2,8)$ ($p=0$, \cite{GAP})}
\begin{tabular}{|c|c|c|c|c|c|c|c|c|c|}
  \hline \hline
           & 1a& 2a & 3a & 7a & 7b & 7c & 9a & 9b & 9c \\ \hline \hline
 $\chi_1$  & 1 & 1  & 1  & 1  & 1  & 1  & 1  & 1  & 1  \\
 $\chi_2$  & 7 & -1 & -2 & 0  & 0  & 0  & 1  & 1  & 1  \\
 $\chi_3$  & 7 & -1 & 1  & 0  & 0  & 0  & $\mathcal{D}$  & $\mathcal{E}$  & $\mathcal{F}$  \\
 $\chi_4$  & 7 & -1 & 1  & 0  & 0  & 0  & $\mathcal{E}$  & $\mathcal{F}$  & $\mathcal{D}$  \\
 $\chi_5$  & 7 & -1 & 1  & 0  & 0  & 0  & $\mathcal{F}$  & $\mathcal{D}$  & $\mathcal{E}$ \\
 $\chi_6$  & 8 & 0 & -1  & 1  & 1  & 1  & -1 & -1 & -1 \\
 $\chi_7$  & 9 & 1 & 0 & $\mathcal{A}$ & $\mathcal{C}$ & $\mathcal{B}$ & 0  &  0 & 0  \\
 $\chi_8$  & 9 & 1 & 0 & $\mathcal{B}$ & $\mathcal{A}$ & $\mathcal{C}$ & 0 & 0 & 0 \\
 $\chi_9$  & 9 & 1 & 0 & $\mathcal{C}$ & $\mathcal{B}$ & $\mathcal{A}$ & 0 & 0 & 0 \\
  \hline
\end{tabular}
\end{center}

where $\mathcal{A}= \alpha^3+\alpha^4$, $\mathcal{B} = \alpha^2+\alpha^5$, $\mathcal{C} = \alpha+\alpha^6$, $\mathcal{D} = -\zeta^4-\zeta^5$, $\mathcal{E} = -\zeta^2-\zeta^7$, $\mathcal{F} = \zeta^2+\zeta^4+\zeta^5+\zeta^7$, $\ds{\alpha=e^{\frac{2 \pi i}{7}}}$ and $\ds{\zeta=e^{\frac{2 \pi i}{9}}}$.
\end{table}

\vspace{-0.25in}

\begin{table}[h!]
  \begin{center}
\caption{Character Table of $PSL(2,8)$ ($p=2$, \cite{GAP})}
\begin{tabular}{|c|c|c|c|c|c|c|c|c|}
  \hline   \hline
      & 1a & 3a & 7a & 7b & 7c & 9a & 9b & 9c\\  \hline  \hline
$\chi_1$   &  1 & 1 & 1 & 1 & 1 & 1 & 1 & 1\\
$\chi_2$   &  2 &-1 & $\mathcal{A}$ & $\mathcal{C}$ & $\mathcal{B}$ & $\mathcal{G}$ & $\mathcal{I}$ & $\mathcal{H}$\\
$\chi_3$   &  2 &-1 & $\mathcal{B}$ & $\mathcal{A}$ & $\mathcal{C}$ & $\mathcal{H}$ & $\mathcal{G}$ & $\mathcal{I}$\\
$\chi_4$   &  2 &-1 & $\mathcal{C}$ & $\mathcal{B}$ & $\mathcal{A}$ & $\mathcal{I}$ & $\mathcal{H}$ & $\mathcal{G}$\\
$\chi_5$   &  4 & 1 & $\mathcal{D}$ & $\mathcal{F}$ & $\mathcal{E}$ & $\mathcal{J}$ & $\mathcal{L}$ & $\mathcal{K}$\\
$\chi_6$   &  4 & 1 & $\mathcal{E}$ & $\mathcal{D}$ & $\mathcal{F}$ & $\mathcal{K}$ & $\mathcal{J}$ & $\mathcal{L}$\\
$\chi_7$   &  4 & 1 & $\mathcal{F}$ & $\mathcal{E}$ & $\mathcal{D}$ & $\mathcal{L}$ & $\mathcal{K}$ & $\mathcal{J}$\\
$\chi_8$   & 8  &-1 & 1 & 1 & 1 &-1 &-1 &-1\\
\hline
\end{tabular}\\

  \end{center}

where $\mathcal{A} = \alpha+\alpha^6$,  $\mathcal{B} = \alpha^3+\alpha^4$,  $\mathcal{C} = \alpha^2+\alpha^5$,  $\mathcal{D} =\alpha+\alpha^3+\alpha^4+\alpha^6$,  $\mathcal{E} = \alpha^2+\alpha^3+\alpha^4+\alpha^5$,  $\mathcal{F} = \alpha+\alpha^2+\alpha^5+\alpha^6$,  $\mathcal{G} = -\zeta^2-\zeta^4-\zeta^5-\zeta^7$,
$\mathcal{H} = \zeta^4+\zeta^5$,  $\mathcal{I} = \zeta^2+\zeta^7$,  $\mathcal{J} = -\zeta^2+\zeta^3-\zeta^4-\zeta^5+\zeta^6-\zeta^7$,
$\mathcal{K} = \zeta^3+\zeta^4+\zeta^5+\zeta^6$, $\mathcal{L}= \zeta^2+\zeta^3+\zeta^6+\zeta^7$, $\ds{\alpha=e^{\frac{2 \pi i}{7}}}$ and $\ds{\zeta=e^{\frac{2 \pi i}{9}}}$

\end{table}

\begin{table}[h!]
  \begin{center}
\caption{Character Table of $PSL(2,8)$  ($p=7$, \cite{GAP})}
\begin{tabular}{|c|c|c|c|c|c|c|}
  \hline   \hline
   & 1a & 2a & 3a & 9a & 9b & 9c\\  \hline  \hline

$\chi_1$ &   1 & 1 & 1 & 1 & 1 & 1\\
$\chi_2$ &   7 &-1 &-2 & 1 & 1 & 1\\
$\chi_3$ &   7 &-1 & 1 & $\mathcal{A}$ & $\mathcal{C}$ & $\mathcal{B}$\\
$\chi_4$ &   7 &-1 & 1 & $\mathcal{B}$ & $\mathcal{A}$ & $\mathcal{C}$\\
$\chi_5$ &   7 &-1 & 1 & $\mathcal{C}$ & $\mathcal{B}$ & $\mathcal{A}$\\
$\chi_6$ &   8 & 0 &-1 &-1 &-1 &-1\\
\hline
\end{tabular}\\[2mm]

  \end{center}

where  $\mathcal{A} = \zeta^2+\zeta^4+\zeta^5+\zeta^7$, $\mathcal{B} = -\zeta^4-\zeta^5$, $\mathcal{C} = -\zeta^2-\zeta^7$ and $\ds{\zeta=e^{\frac{2 \pi i}{9}}}$

\end{table}

\section{Proof of Theorem 2}

Let $G=PSL(2,17)$. Clearly $|G|=2448 = 2^{4} \cdot 3^{2} \cdot 17$ and $exp(G) = 1224 = 2^{3} \cdot 3^{2} \cdot 17$.  Initially for any torsion unit
of $V(\ZZ G)$ of order $k$, we have that
 \[\nu_{2a} + \nu_{3a} +\nu_{4a} +\nu_{8a} +\nu_{8b} +\nu_{9a} +\nu_{9b} +\nu_{9c} +\nu_{17a} +\nu_{17b}=1.\]

By Proposition \ref{CL}, we need only to consider torsion units of $V(\ZZ G)$ of order $2$, $3$, $4$, $6$, $8$, $9$, $17$, $34$ and $51$.\\

Case $(i)$. Let $u$ be a torsion unit in $V(\ZZ G)$ of order $2$ and $3$.  Clearly $u$ is rationally conjugated to some element $g \in G$ by Propositions \ref{P:1} and \ref{P:3}.\\

Case $(ii)$. Let $u$ be a torsion unit in $V(\ZZ G)$ of order $17$.  Using Proposition \ref{6:1}, $u$ is rationally conjugated to some element $g \in G$.\\

Case $(iii)$. Let $u$ be a torsion unit in $V(\ZZ G)$ of order $4$.  Using Propositions \ref{P:1} \& \ref{P:2},
$\nu_{2a}+\nu_{4a}=1$.  Applying Proposition \ref{P:4} to the character table (Table 5) , we obtain

\[
\begin{split}
\mu_{2}(u,\chi_{2},17) & = \textstyle \frac{1}{4} (2 \gamma_1 + 2) \geq 0; \quad
\mu_{0}(u,\chi_{2},17)  = \textstyle \frac{1}{4} (-2 \gamma_1 + 2) \geq 0,
\end{split}
\]

\noindent where $\gamma_1 =  \nu_{2a} -  \nu_{4a}$.  Clearly $\gamma_1 \in \{ -1, 1 \}$.  It follows that
the only possible integer values for $(\nu_{2a},\nu_{4a})$ are $(1,0)$ and $(0,1)$. Therefore $u$ is rationally conjugated to some element $g \in G$ by Proposition \ref{P:3}.\\

Case $(iv)$. Let $u$ be a torsion unit in $V(\ZZ G)$ of order $8$.  Using Propositions \ref{P:1} \& \ref{P:2},
\[
\nu_{2a}+\nu_{4a}+\nu_{8a}+\nu_{8b}=1.
\]

\noindent Applying Proposition \ref{P:4} to the character tables (Table 4 \& 5) , we obtain

\[
\begin{split}
\mu_{4}(u,\chi_{4},17) & = \textstyle \frac{1}{8} (4 \gamma_1 + 4) \geq 0; \quad
\mu_{0}(u,\chi_{4},17)  = \textstyle \frac{1}{8} (-4 \gamma_1 + 4) \geq 0; \\
\mu_{0}(u,\chi_{9},*) & = \textstyle \frac{1}{8} (8 \gamma_2 +16) \geq 0; \quad
\mu_{4}(u,\chi_{9},*)  = \textstyle \frac{1}{8} (-8 \gamma_2 +16) \geq 0; \\
\mu_{1}(u,\chi_{2},17) & = \textstyle \frac{1}{8} (4 \gamma_3 + 4) \geq 0; \quad
\mu_{3}(u,\chi_{2},17)  = \textstyle \frac{1}{8} (- 4\gamma_3+ 4) \geq 0; \\
\mu_{4}(u,\chi_{2},17)  &= \textstyle \frac{1}{8} (4 \gamma_4 + 4) \geq 0; \quad
\mu_{4}(u,\chi_{3},17)  = \textstyle \frac{1}{8} (4\gamma_5 + 4) \geq 0,
\end{split}
\]
where  $\gamma_1 =  \nu_{2a} +  \nu_{4a} -  \nu_{8a} -  \nu_{8b}$, $\gamma_2 =  \nu_{2a} -  \nu_{4a}$
and $\gamma_3 = \nu_{8a} -  \nu_{8b}$, $\gamma_4 = \nu_{2a} - \nu_{4a} - \nu_{8a} -  \nu_{8b}$ and $\gamma_5=
 -\nu_{2a} +  \nu_{4a} - \nu_{8a} -  \nu_{8b}$.  Now $\gamma_1 \in \{ -1, 1 \}$, $\gamma_2 \in \{ -2, -1, 0, 1, 2 \}$, $\gamma_3 \in \{ -1, 1 \}$.   It follows that the only possible integer values for $(\nu_{2a},\nu_{4a},\nu_{8a},\nu_{8b})$ are $(0, 0, 0, 1)$ and $(0, 0, 1, 0)$. Therefore $u$ is rationally conjugated to some element $g \in G$ by Proposition \ref{P:3}.\\

Case $(v)$. Let $|u|=9$.  Using Propositions \ref{P:1} \& \ref{P:2},
\[
\nu_{3a}+\nu_{9a}+\nu_{9b}+\nu_{9c}=1.
\]
Applying Proposition \ref{P:4} to the character tables (Table 4 \& 5) , we obtain
\[
\begin{split}
\mu_{3}(u,\chi_{4},*) & = \textstyle \frac{1}{9} (3\gamma_1 + 12) \geq 0; \quad
\mu_{0}(u,\chi_{4},*)  = \textstyle \frac{1}{9} (-6\gamma_1 + 12) \geq 0; \\
\mu_{2}(u,\chi_{2},17) & = \textstyle \frac{1}{9} (3\gamma_2 + 3) \geq 0; \quad
\mu_{1}(u,\chi_{3},17)  = \textstyle \frac{1}{9} (- 3\gamma_2 + 6) \geq 0; \\
\mu_{2}(u,\chi_{3},17) & = \textstyle \frac{1}{9} (3 \gamma_3 + 6) \geq 0; \quad
\mu_{4}(u,\chi_{2},17)  = \textstyle \frac{1}{9} (- 3\gamma_3 + 3) \geq 0; \\
\mu_{1}(u,&\chi_{2},17)  = \textstyle \frac{1}{9} (- 3\gamma_4 + 3) \geq 0,
\end{split}
\]
where $\gamma_1 = 2 \nu_{3a} -  \nu_{9a} -  \nu_{9b} -  \nu_{9c}$, $\gamma_2 =  2 \nu_{9a} -  \nu_{9b} -  \nu_{9c}$, $\gamma_3 =   \nu_{9a} - 2 \nu_{9b} +  \nu_{9c}$ and $\gamma_4=\nu_{9a} +\nu_{9b} -2 \nu_{9c}$.  Now  $\gamma_1 \in \{ -4, -1, 2 \}$, $\gamma_2 \in \{ -1, 2 \}$ and $\gamma_3 \in \{ -2, 1 \}$.  It follows that the only possible integer values for $(\nu_{3a},\nu_{9a},\nu_{9b},\nu_{9c})$ are $( 0, 0, 0, 1 )$, $(0, 0, 1, 0)$ and $(0, 1, 0, 0)$. Therefore $u$ is rationally conjugated to some element $g \in G$ by Proposition \ref{P:3}.\\

Case $(vi)$. Let us consider all possible values of $|u|$ which do not appear in $G$.  By \cite{CL},
$|u| \in \{6,34,51\}$. By proposition \ref{6:6}, there doesn't exist any torsion units of order $6$ (since $G$
doesn't contain any such elements).  Finally if $|u| \in \{34,51\}$, then $u$ is rationally conjugated to some element $g \in G$ by Proposition \ref{6:3}. This completes the proof.

\begin{table}[h!]
  \begin{center}
 \caption{Character Table of $PSL(2,17)$ ($p=0$, \cite{GAP})}
\begin{tabular}{|c|c|c|c|c|c|c|c|c|c|c|c|}
  \hline  \hline
                & 1a & 2a & 3a & 4a& 8a & 8b & 9a & 9b & 9c & 17a & 17b \\   \hline \hline
$\chi_1$        & 1  & 1  & 1  & 1 & 1  &  1 &  1 &  1 &  1 &  1  & 1 \\
$\chi_2$        & 9  & 1  & 0  & 1 &-1  & -1 &  0 &  0 &  0 &  $\mathcal{E}$  & $\mathcal{F}$ \\
$\chi_3$        & 9  & 1  & 0  & 1 &-1  & -1 &  0 &  0 &  0 &  $\mathcal{F}$  & $\mathcal{E}$ \\
$\chi_4$        & 16 & 0  & -2 & 0 & 0  &  0 &  1 &  1 &  1 & -1  & -1 \\
$\chi_5$        & 16 & 0  & 1  & 0 & 0  &  0 &  $\mathcal{B}$ &  $\mathcal{C}$ &  $\mathcal{D}$ & -1  & -1 \\
$\chi_6$        & 16 & 0  & 1  & 0 & 0  &  0 &  $\mathcal{C}$ &  $\mathcal{D}$ &  $\mathcal{B}$ & -1  & -1 \\
$\chi_7$        & 16 & 0  & 1  & 0 & 0  &  0 &  $\mathcal{D}$ &  $\mathcal{B}$ &  $\mathcal{C}$ & -1  & -1 \\
$\chi_8$        & 17 & 1  & -1 & 1 & 1  &  1 & -1 &  -1& -1 &  0  & 0 \\
$\chi_9$        & 18 & 2  & 0  & -2& 0  &  0 &  0 &  0 &  0 &  1  & 1 \\
$\chi_{10}$     & 18 & -2 & 0  & 0 & $\mathcal{A}$  & -$\mathcal{A}$ &  0 &  0 &  0 &  1  & 1 \\
$\chi_{11}$     & 18 & -2 & 0  & 0 & -$\mathcal{A}$ &  $\mathcal{A}$ &  0 &  0 &  0 &  1  & 1 \\
  \hline
\end{tabular}
\\[2mm]
\end{center}
where  $\mathcal{A} = -\alpha+\alpha^3$, $\mathcal{B} = -\zeta^2-\zeta^7$, $\mathcal{C} = -\zeta^4-\zeta^5$, $\mathcal{D} = \zeta^2+\zeta^4+\zeta^5+\zeta^7$,
$\mathcal{E} = -\delta-\delta^2-\delta^4-\delta^8-\delta^9-\delta^{13}-\delta^{15}-\delta^{16}$, $\mathcal{F} = -\delta^3-\delta^5-\delta^6-\delta^7-\delta^{10}-\delta^{11}-\delta^{12}-\delta^{14}$,
$\ds{\alpha=e^{\frac{2 \pi i}{8}}}$,  $\ds{\zeta=e^{\frac{2 \pi i}{9}}}$ and $\ds{\delta=e^{\frac{2 \pi i}{17}}}$.
\end{table}

\begin{table}[h!]
  \begin{center}
 \caption{Character Table of $PSL(2,17)$  ($p=17$, \cite{GAP})}
\begin{tabular}{|c|c|c|c|c|c|c|c|c|c|}
  \hline \hline
     & 1a & 2a & 3a & 4a & 8a & 8b & 9a & 9b & 9c\\ \hline \hline
$\chi_1$  &  1 & 1  & 1 & 1  &1  & 1 & 1 & 1 & 1\\
$\chi_2$  &  3 & -1 & 0 & 1  & $\mathcal{A}$ & $\mathcal{K}$ & $\mathcal{B}$ & $\mathcal{E}$ & $\mathcal{H}$\\
$\chi_3$  &  5 & 1  &-1 & -1 & $\mathcal{A}$ & $\mathcal{K}$ & $\mathcal{C}$ & $\mathcal{F}$ & $\mathcal{I}$\\
$\chi_4$  &  7 &-1  & 1 & -1 & 1 &  1 & $\mathcal{D}$ & $\mathcal{G}$ & $\mathcal{J}$\\
$\chi_5$  &  9 & 1  & 0 & 1  &-1 & -1 & 0 & 0 & 0\\
$\chi_6$  & 11 &-1  & -1& 1  &-$\mathcal{A}$ & -$\mathcal{K}$ & -$\mathcal{D}$ &-$\mathcal{G}$ &-$\mathcal{J}$\\
$\chi_7$  & 13 & 1  & 1 &-1  &-$\mathcal{A}$ & -$\mathcal{K}$ & -$\mathcal{C}$ &-$\mathcal{F}$ &-$\mathcal{I}$\\
$\chi_8$  & 15 &-1  & 0 &-1  &-1 & -1 & -$\mathcal{B}$ &-$\mathcal{E}$ &-$\mathcal{H}$\\
$\chi_9$  & 17 & 1  &-1 & 1  &1  & 1 & -1 &-1 &-1\\
\hline
\end{tabular}
\\[2mm]
\end{center}
where $\mathcal{A} = 1+\alpha-\alpha^3$, $\mathcal{B} = \zeta^2-\zeta^3-\zeta^6+\zeta^7$, $\mathcal{C} = \zeta^2-\zeta^3+\zeta^4+\zeta^5-\zeta^6+\zeta^7$,
$\mathcal{D} = \zeta^2+\zeta^4+\zeta^5+\zeta^7$, $\mathcal{E} = -\zeta^3+\zeta^4+\zeta^5-\zeta^6$, $\mathcal{F} = -\zeta^2-\zeta^3-\zeta^6-\zeta^7$,
$\mathcal{G} = -\zeta^2-\zeta^7$, $\mathcal{H} = -\zeta^2-\zeta^3-\zeta^4-\zeta^5-\zeta^6-\zeta^7$, $\mathcal{I} = -\zeta^3-\zeta^4-\zeta^5-\zeta^6$,
$\mathcal{J} = -\zeta^4-\zeta^5$, $\mathcal{K} =1-\alpha+\alpha^3$,  $\ds{\alpha=e^{\frac{2 \pi i}{8}}}$ and $\ds{\zeta=e^{\frac{2 \pi i}{9}}}$.
\end{table}


\begin{bibdiv}
  \begin{biblist}


\bib{AB}{article}{
   author={Artamonov, V. A.},
   author={Bovdi, A. A.},
   title={Integral group rings: groups of invertible elements and classical
   $K$-theory},
   language={Russian},
   note={Translated in J. Soviet Math.\ {\bf 57} (1991), no.\ 2,
   2931--2958},
   conference={
      title={Algebra. Topology. Geometry, Vol.\ 27 (Russian)},
   },
   book={
      series={Itogi Nauki i Tekhniki},
      publisher={Akad. Nauk SSSR Vsesoyuz. Inst. Nauchn. i Tekhn. Inform.},
      place={Moscow},
   },
   date={1989},
   pages={3--43, 232},
   review={},
}




\bib{BGK}{article}{
   author={Bovdi, V.},
   author={Grishkov, A.},
   author={Konovalov, A.},
   title={Kimmerle conjecture for the Held and O'Nan sporadic simple groups},
   journal={Sci. Math. Jpn.},
   volume={69},
   date={2009},
   number={3},
   pages={353--361},
   issn={},
   review={},
}

\bib{BH}{article}{
   author={Bovdi, V. },
   author={Hertweck, M.},
   title={Zassenhaus conjecture for central extensions of $S\sb 5$},
   journal={J. Group Theory},
   volume={11},
   date={2008},
   number={1},
   pages={63--74},
   issn={1433-5883},
   review={},
   doi={},
}


\bib{BHK}{article}{
   author={Bovdi, V.},
   author={H{\"o}fert, C.},
   author={Kimmerle, W.},
   title={On the first Zassenhaus conjecture for integral group rings},
   journal={Publ. Math. Debrecen},
   volume={65},
   date={2004},
   number={3-4},
   pages={291--303},
   issn={},
   review={},
}

\bib{BJK}{article}{
   author={Bovdi, V.},
   author={Jespers, E.},
   author={Konovalov, A. B.},
   title={Torsion units in integral group rings of Janko simple groups},
   journal={Math. Comp.},
   volume={80},
   date={2011},
   number={273},
   pages={593--615},
   issn={0025-5718},
   review={\MR{2728996 (2011j:20010)}},
   doi={10.1090/S0025-5718-2010-02376-2},
}


\bib{BKM11}{article}{
   author={Bovdi, V.},
   author={Konovalov, A. B.},
   title={Integral group ring of the first Mathieu simple group},
   conference={
      title={Groups St. Andrews 2005. Vol. 1},
   },
   book={
      series={London Math. Soc. Lecture Note Ser.},
      volume={339},
      publisher={Cambridge Univ. Press},
      place={Cambridge},
   },
   date={2007},
   pages={237--245},
   review={},
   doi={},
}


\bib{BKM23}{article}{
   author={ Bovdi, V.},
   author={Konovalov, A. B.},
   title={Integral group ring of the Mathieu simple group $M\sb {23}$},
   journal={Comm. Algebra},
   volume={36},
   date={2008},
   number={7},
   pages={2670--2680},
   issn={},
   review={},
   doi={},
}


\bib{BKR}{article}{
   author={Bovdi, V.},
   author={Konovalov, A. B.},
   title={Integral group ring of Rudvalis simple group},
   language={English, with English and Ukrainian summaries},
   journal={Ukra\"\i n. Mat. Zh.},
   volume={61},
   date={2009},
   number={1},
   pages={3--13},
   issn={1027-3190},
   translation={
      journal={Ukrainian Math. J.},
      volume={61},
      date={2009},
      number={1},
      pages={1--13},
      issn={},
  review={},
   doi={},
}
}

\bib{BKHS}{article}{
   author={ Bovdi, V.},
   author={ Konovalov, A. B.},
   title={Torsion units in integral group ring of Higman-Sims simple group},
   journal={Studia Sci. Math. Hungar.},
   volume={47},
   date={2010},
   number={1},
   pages={1--11},
   issn={},
   review={},
   doi={},
}


\bib{BKML}{article}{
   author={Bovdi, V.},
   author={Konovalov, A. B.},
   title={Integral group ring of the McLaughlin simple group},
   journal={Algebra Discrete Math.},
   date={2007},
   number={2},
   pages={43--53},
   issn={},
   review={},
}


\bib{BKM24}{article}{
   author={  Bovdi, V.},
   author={Konovalov, A. B.},
   title={Integral group ring of the Mathieu simple group $M\sb {24}$},
   journal={J. Algebra Appl.},
   volume={11},
   date={2012},
   number={1},
   pages={1250016, 10},
   issn={},
   review={},
   doi={},
}

\bib{BKL}{article}{
   author={Bovdi, V.},
   author={Konovalov, A. B.},
   author={Linton, S.},
   title={Torsion units in integral group ring of the Mathieu simple group
   ${\rm M}\sb {22}$},
   journal={LMS J. Comput. Math.},
   volume={11},
   date={2008},
   pages={28--39},
   issn={},
   review={},
}

\bib{BKLCW}{article}{
   author={Bovdi, V. },
   author={Konovalov, A. B.},
   author={Linton, S.},
   title={Torsion units in integral group rings of Conway simple groups},
   journal={Internat. J. Algebra Comput.},
   volume={21},
   date={2011},
   number={4},
   pages={615--634},
   issn={},
   review={},
   doi={},
}


\bib{BKM}{article}{
   author={Bovdi, V. A.},
   author={Konovalov, A. B.},
   author={Marcos, E. N.},
   title={Integral group ring of the Suzuki sporadic simple group},
   journal={Publ. Math. Debrecen},
   volume={72},
   date={2008},
   number={3-4},
   pages={487--503},
   issn={},
   review={},
}

\bib{BKS}{article}{
   author={Bovdi, V.},
   author={Konovalov, A. B.},
   author={Siciliano, S.},
   title={Integral group ring of the Mathieu simple group $M\sb {12}$},
   journal={Rend. Circ. Mat. Palermo (2)},
   volume={56},
   date={2007},
   number={1},
   pages={125--136},
   issn={},
   review={},
   doi={},
}




\bib{CL}{article}{
   author={Cohn, J. A.},
   author={Livingstone, D.},
   title={On the structure of group algebras. I},
   journal={Canad. J. Math.},
   volume={17},
   date={1965},
   pages={583--593},
   issn={},
   review={},
}

\bib{MH1}{article}{
   author={Hertweck, M.},
   title={On the torsion units of some integral group rings},
   journal={Algebra Colloq.},
   volume={13},
   date={2006},
   number={2},
   pages={329--348},
   issn={},
   review={},
}

\bib{MH2}{article}{
   author={ Hertweck, M.},
   title={Zassenhaus conjecture for $A\sb 6$},
   journal={Proc. Indian Acad. Sci. Math. Sci.},
   volume={118},
   date={2008},
   number={2},
   pages={189--195},
   issn={},
   review={},
   doi={},
}



\bib{MH3}{article}{
   author={Hertweck, M.},
   title={Partial augmentations and Brauer character values of torsion units in group rings},
   journal={(http://arxiv.org/abs/math/0612429)},
   volume={},
   date={2007},
   number={},
   pages={},
   issn={},
   review={},
   doi={},
}

\bib{MH4}{article}{
   author={  Hertweck, M.},
   title={Contributions to the integral representation theory of groups (Habilitationsschrift, University of Stuttgart)},
   journal={electronic publication, http://elib.uni-stuttgart.de/opus/volltexte/2004/1638.},
   volume={},
   date={2004},
   number={},
   pages={191 pages},
   issn={},
   review={},
   doi={},
}


\bib{HHK}{article}{
   author={Hertweck, M.},
   author={H{\"o}fert, C.},
   author={Kimmerle, W.},
   title={Finite groups of units and their composition factors in the
   integral group rings of the group ${\rm PSL}(2,q)$},
   journal={J. Group Theory},
   volume={12},
   date={2009},
   number={6},
   pages={873--882},
   issn={ 1433-5883},
}

\bib{HK}{article}{
   author={H{\"o}fert, C.},
   author={Kimmerle, W.},
   title={On torsion units of integral group rings of groups of small order},
   conference={
      title={Groups, rings and group rings},
   },
   book={
      series={Lect. Notes Pure Appl. Math.},
      volume={248},
      publisher={Chapman \& Hall/CRC, Boca Raton, FL},
   },
   date={2006},
   pages={243--252},
}





\bib{GAP}{article}{
   author={The GAP Group},
   title={GAP- Groups, Algorithms and Programming},
   journal={version 4.4},
   volume={2006 (http:/www.gap-system.org)},
   date={},
   number={},
   pages={},
   issn={},
   review={},
   doi={},
}

\bib{WK}{article}{
   author={Kimmerle, W.},
   title={On the prime graph of the unit group of integral group rings of
   finite groups},
   conference={
      title={Groups, rings and algebras},
   },
   book={
      series={Contemp. Math.},
      volume={420},
      publisher={Amer. Math. Soc.},
      place={Providence, RI},
   },
   date={2006},
   pages={215--228},
   review={\MR{2279241 (2007m:16040)}},
}


\bib{lag}{article}{
   author={Bovdi, A.},
   author={Konovalov, A.},
   author={R. Rossmanith},
   author={C. Schneider},
   title={LAGUNA --- Lie AlGebras and UNits of group Algebras},
   journal={Version 3.5.0; 2009},
   volume={},
   date={http://www.cs.st-andrews.ac.uk/~alexk/laguna.htm},
   number={},
   pages={},
   issn={},
   review={},
   doi={},
}



\bib{LP1}{article}{
   author={Luthar, I. S.},
   author={Passi, I. B. S.},
   title={Zassenhaus conjecture for $A\sb 5$},
   journal={Proc. Indian Acad. Sci. Math. Sci.},
   volume={99},
   date={1989},
   number={1},
   pages={1--5},
   issn={},
   review={},
   doi={},
}

\bib{LP2}{article}{
   author={Luthar, I. S.},
   author={Trama, Poonam},
   title={Zassenhaus conjecture for $S\sb 5$},
   journal={Comm. Algebra},
   volume={19},
   date={1991},
   number={8},
   pages={2353--2362},
   issn={},
   review={},
   doi={},
}

\bib{MRSW}{article}{
   author={Marciniak, Z.},
   author={Ritter, J.},
   author={Sehgal, S. K.},
   author={Weiss, A.},
   title={Torsion units in integral group rings of some metabelian groups.
   II},
   journal={J. Number Theory},
   volume={25},
   date={1987},
   number={3},
   pages={340--352},
   issn={},
   review={},
   doi={},
}



\bib{OB}{article}{
   title={Mini-Workshop},
   journal={Arithmetik von Gruppenringen, Oberwolfach Rep. 4 (2007), no. 4, 3209 $–-$
3239, Abstracts from the mini-workshop held November 25 $-$ December 1, 2007, Organized by
E. Jespers, Z. Marciniak, G. Nebe and W. Kimmerle, Oberwolfach Reports. Vol. 4, no. 4.},
   publisher={}
   volume={},
   date={},
   pages={},
   review={\MR{MR2463649}},
}


\bib{RS}{article}{
   author={Roggenkamp, K.},
   author={Scott, L.},
   title={Isomorphisms of $p$-adic group rings},
   journal={Ann. of Math. (2)},
   volume={126},
   date={1987},
   number={3},
   pages={593--647},
   issn={0003-486X},
   review={\MR{916720 (89b:20021)}},
   doi={10.2307/1971362},
}

\bib{MS}{article}{
   author={Salim, M.},
   title={Torsion units in the integral group ring of the alternating group
   of degree 6},
   journal={Comm. Algebra},
   volume={35},
   date={2007},
   number={12},
   pages={4198--4204},
   issn={0092-7872},
   doi={10.1080/00927870701545069},
}



\bib{W}{article}{
   author={Weiss, A.},
   title={Rigidity of $p$-adic $p$-torsion},
   journal={Ann. of Math. (2)},
   volume={127},
   date={1988},
   number={2},
   pages={317--332},
   issn={0003-486X},
   review={\MR{932300 (89g:20010)}},
   doi={10.2307/2007056},
}

\bib{HZ}{article}{
   author={Zassenhaus, H.},
   title={On the torsion units of finite group rings},
   conference={
      title={Studies in mathematics (in honor of A. Almeida Costa)
      (Portuguese)},
   },
   book={
      publisher={Instituto de Alta Cultura},
      place={Lisbon},
   },
   date={1974},
   pages={119--126},
   review={},
}


  \end{biblist}
\end{bibdiv}

\end{document}